\documentclass[11pt]{amsart}
\usepackage{graphicx}
\usepackage{amssymb}
\begin{document}

\title{Minkowski Bisectors, Minkowski Cells, and Lattice Coverings}

\author{Chuanming Zong}

\thanks{{\em 2010 Mathematics Subject Classification}: primary 52C17; secondary 52B10, 52C07}

\medskip
\bigskip

\begin{abstract}
As a counterpart of Hilbert's 18th problem, it is natural to raise the question to determine the density of the thinnest covering of $\mathbb{E}^3$ by identical copies of a given geometric object, such as a unit ball, a regular tetrahedron or a regular octahedron. However, systematic study on this problem and its generalizations was started much later comparing with that on packings and the known knowledge on coverings is very limited. By studying Minkowski bisectors and Minkowski cells, this paper introduces a new way to study the density $\theta^*(C)$ of the thinnest lattice covering of $\mathbb{E}^n$ by a centrally symmetric convex body $C$.  Some basic results and unexpected phenomena such as Example 1, Theorems 2 and 4, and Corollary 1 about Minkowski bisectors, Minkowski cells and covering densities are discovered. Three basic problems about Minkowski cells and parallelohedra are presented.
\end{abstract}

\maketitle

\medskip
\bigskip
\section{Minkowski Metrics and Centrally Symmetric Convex Bodies}

\bigskip
Let $\mathbb{E}^n$ denote the $n$-dimensional Euclidean space with an orthonormal basis $\{ {\bf e}_1, {\bf
e}_2, \cdots , {\bf e}_n\}$. For convenience we use small letters to denote real numbers, use small bold letters to denote points (or vectors) and use capital letters to denote sets of points in the space. In particular,
let ${\bf o}$ denote the origin of $\mathbb{E}^n$, and let $B_n$ denote the $n$-dimensional unit ball $\{{\bf x}:\ \sum |x_i|^2\le 1\}$.

A metric $\| \cdot \|$ define in $\mathbb{E}^n$ is called a {\it Minkowski
metric} if it satisfies the following conditions.

\medskip
\noindent {\bf 1.} {\it $0\le \| {\bf x}\|<\infty$ holds
for all vectors ${\bf x}\in \mathbb{E}^n$ and the equality holds if and
only if ${\bf x}={\bf o}.$}

\noindent {\bf 2.} {\it $\| \lambda {\bf x}\|=|\lambda|\cdot \|
{\bf x}\|$ holds for all vectors ${\bf x}\in \mathbb{E}^n$ and all real
numbers $\lambda $.}

\noindent {\bf 3.} {\it $\| {\bf x}+{\bf y}\|\le \|{\bf x}\|+\|
{\bf y}\|$ holds for all vectors ${\bf x},$ ${\bf y}\in \mathbb{E}^n$.}

\medskip
Let us define
$$C=\{ {\bf x}\in \mathbb{E}^n:\ \| {\bf x}\|\le 1\}.$$
Since $\| -{\bf x}\|=\| {\bf x}\|$, if ${\bf x}\in C$, then $-{\bf
x}\in C$. Therefore $C$ is centrally symmetric. On the other hand,
if both ${\bf x}$ and ${\bf y}$ belong to $C$ and if $\mu$ is
a real number satisfying $0<\mu <1$, by the third condition of the
metric we have
$$\| \mu {\bf x}+(1-\mu ){\bf y}\|\le \mu \| {\bf
x}\|+(1-\mu )\| {\bf y}\| \le \mu +(1-\mu )=1.$$
In other words, $C$ is convex. In addition, the origin ${\bf o}$ is an
interior point of $C$. Otherwise there would be a point ${\bf x}$
satisfying $\| {\bf x}\|=\infty $. Finally, since $\| {\bf x}\|=0$
holds if and only if ${\bf x}={\bf o}$, it can be deduced that $C$
is bounded. As a conclusion, $C$ is a centrally symmetric convex
body centered at the origin ${\bf o}$.

On the other hand, each centrally symmetric convex body $C$ centered at the origin ${\bf o}$ can define
a Minkowski metric. Let $\overline{{\bf x}}$ denote the point in the direction of ${\bf x}$ and on the boundary of $C$, and define
$$\| {\bf x}\|=\left\{
\begin{array}{cl}
0 &\mbox{if ${\bf
x}= {\bf o}$,}\\
\lambda &\mbox{where ${\bf x}=\lambda \overline{\bf x}$}.
\end{array}\right.$$
In other words, $\| {\bf x}\|$ is the positive number $\lambda$ that ${\bf x}$ is on the boundary of $\lambda C$.

It can be easily verified that
$$0<\| {\bf x}\| <\infty $$
holds for all ${\bf x}\not= {\bf o}$, and
$$\|\lambda {\bf x}\|=|\lambda |\cdot \| {\bf x}\|$$
holds for all vectors ${\bf x}\in \mathbb{E}^n$ and all real numbers $\lambda $. Since $C$ is convex, for any pair of vectors ${\bf x}$ and
${\bf y}$ it follows that
$${{{\bf x}+{\bf y}}\over {\| {\bf x}\|+\|{\bf y}\|}}={{\| {\bf x}\|}\over
{\| {\bf x}\|+\|{\bf y}\|}}\overline{\bf x} +{{\|{\bf y}\|}\over
{\| {\bf x}\|+\|{\bf y}\|}}\overline{\bf y}\in C,$$

$$\left\|{{{\bf x}+{\bf y}}\over {\| {\bf x}\|+\|{\bf
y}\|}}\right\|={{\|{\bf x}+{\bf y}\|}\over {\| {\bf x}\|+\|{\bf
y}\|}}\le 1$$ and therefore
$$\|{\bf x}+{\bf y}\|\le \| {\bf x}\|+\|{\bf y}\|.$$
Thus $\| \cdot \|$ defines a Minkowski metric on $\mathbb{E}^n$.
In particular, when $C$ is the $n$-dimensional unit ball, the
corresponding Minkowski metric is the Euclidean metric.

As a conclusion of this section, we get the following well-known assertion (see page 7 of \cite{gru87}).

\medskip
\noindent
{\bf Theorem 0.} {\it There is a one-to-one correspondence
between Minkowski metrics and centrally symmetric convex bodies in
$\mathbb{E}^n$. }

\medskip
\section{The Minkowski Bisectors}

\bigskip
Let $\| \cdot \|_C$ denote a Minkowski metric determined by a
centrally symmetric convex body $C$ and let $\| {\bf x}, {\bf
y}\|_C=\| {\bf y}-{\bf x}\|_C$ denote the {\it Minkowski distance}
between two points ${\bf x}$ and ${\bf y}$ with respect to the
metric. For two distinct points ${\bf p}$ and ${\bf q}$ in $\mathbb{E}^n$
we define
$$H({\bf p}, {\bf q})=\left\{ {\bf y}\in \mathbb{E}^n:\ \langle {\bf y},
{\bf q}-{\bf p}\rangle =0\right\}$$ and
$$L({\bf p}, {\bf q}, {\bf x})=\left\{ {\bf x}+\lambda ({\bf q}-{\bf p}):\
\lambda \in \mathbb{R}\right\},$$ where ${\bf x}\in H({\bf p}, {\bf q})$.
In other words, $H({\bf p}, {\bf q})$ is a hyperplane containing the origin
and perpendicular to the vector ${\bf q}-{\bf p}$, and $L({\bf p}, {\bf q}, {\bf x})$
is the straight line determined by the point ${\bf x}$ and the vector ${\bf q}-{\bf
p}$. Then, for every ${\bf x}\in H({\bf p}, {\bf q})$, we define
$$S(C, {\bf p}, {\bf q}, {\bf x})=\left\{ {\bf y}\in L({\bf p}, {\bf q}, {\bf
x}): \ \| {\bf p}, {\bf y}\|_C=\| {\bf q}, {\bf y}\|_C\right\}.$$

Now let us introduce a basic result about $S(C, {\bf p}, {\bf q}, {\bf
x})$ which is essential for the definition of the Minkowski bisectors.

\medskip\noindent
{\bf Lemma 1 (Horv\'ath \cite{hor00}).} {\it For every ${\bf x}\in
H({\bf p}, {\bf q})$, the set $S(C, {\bf p}, {\bf q}, {\bf x})$ is
either a single point or a closed segment.}

\medskip\noindent
{\bf Proof.} Since $C$ is a centrally symmetric convex body, both
$f({\bf y})=\| {\bf p}, {\bf y}\|_C$ and $g({\bf y})=\| {\bf q},
{\bf y}\|_C$ are continuous functions of ${\bf y}$. On the other
hand, when $\lambda $ is sufficiently large, we have
$$f({\bf x}+\lambda ({\bf q}-{\bf p}))>g({\bf x}+\lambda ({\bf q}-{\bf
p}))$$ and
$$f({\bf x}-\lambda ({\bf q}-{\bf p}))<g({\bf x}-\lambda ({\bf q}-{\bf
p})).$$ Thus it follows that $S(C, {\bf p}, {\bf q}, {\bf x})$ is
nonempty. In addition, it can be easily deduced that $S(C, {\bf
p}, {\bf q}, {\bf x})$ is always a compact set.

Now we show that if both ${\bf u}$ and ${\bf v}$ belong to $S(C,
{\bf p}, {\bf q}, {\bf x})$ (see Figure 1) then the whole segment
$[{\bf u}, {\bf v}]$ belongs to $S(C, {\bf p}, {\bf q}, {\bf x})$.

\begin{figure}[ht]
\centering
\includegraphics[height=4.1cm,width=8cm,angle=0]{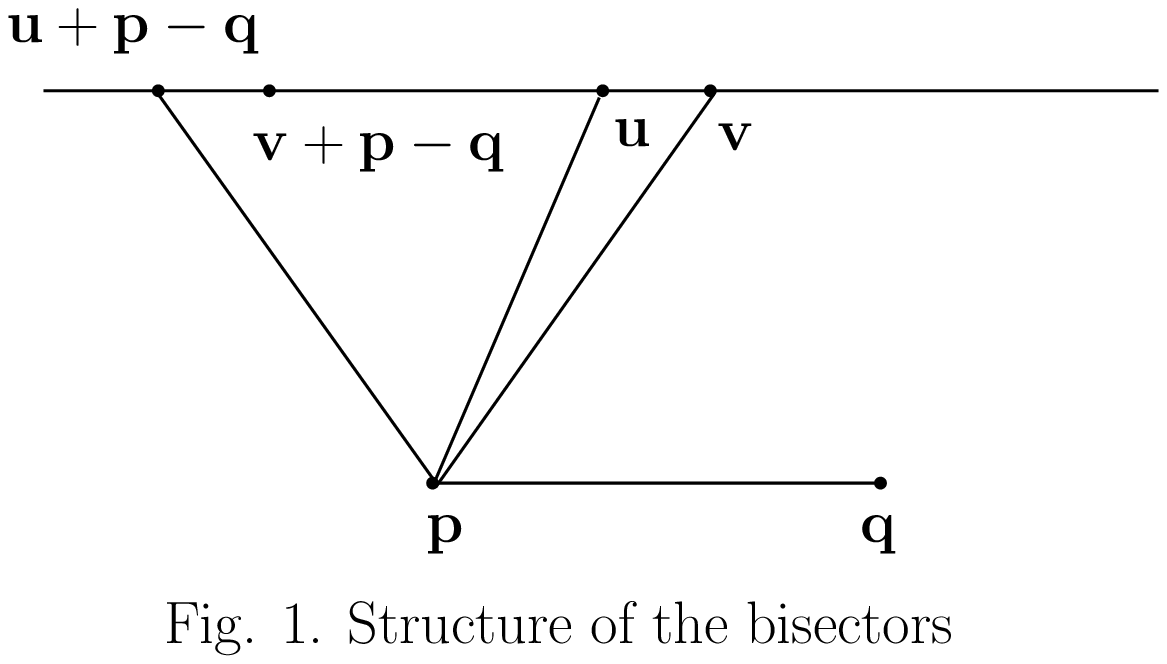}
\end{figure}

\medskip
Assume that
$$\| {\bf p}-{\bf u}\|_C=\| {\bf q}-{\bf u}\|_C=\alpha ,$$
$$\| {\bf p}-{\bf v}\|_C=\| {\bf q}-{\bf v}\|_C=\beta $$
and $\alpha <\beta $. It is easy to see that both ${\bf u}+{\bf
p}-{\bf q}$ and ${\bf v}+{\bf p}-{\bf q}$ belong to $L({\bf p},
{\bf q}, {\bf x})$, ${\bf u}$ is between ${\bf v}$ and ${\bf u}+{\bf p}-{\bf q}$ and therefore
$${\bf u}-{\bf p}=\mu ({\bf v}-{\bf p})+(1-\mu )({\bf
u}+{\bf p}-{\bf q}-{\bf p})$$ holds for some $\mu $ with
$0<\mu <1$. Thus we get
\begin{eqnarray*}
\alpha\!\!\! &=&\!\!\! \| {\bf p}-{\bf u}\|_C\ge \mu \| {\bf
p}-{\bf v}\|_C +(1-\mu )\| {\bf p}-{\bf u}\|_C\\
\!\!\! &\ge & \!\!\! \mu \beta +(1-\mu )\alpha =\alpha
+\mu (\beta -\alpha )>\alpha.
\end{eqnarray*}
By this contradiction we can conclude that $\alpha =\beta $. Then
it follows by convexity (since all ${\bf u}+{\bf p}-{\bf q}$,
${\bf v}+{\bf p}-{\bf q}$, ${\bf u}$ and ${\bf v}$ belong to the
boundary of ${\bf p}+\alpha C$) that the whole segment $[{\bf
u}+{\bf p}-{\bf q}, {\bf v}]$ belongs to the boundary of ${\bf
p}+\alpha C$, the whole segment $[{\bf u}, {\bf v}+{\bf q}-{\bf
p}]$ belongs to the boundary of ${\bf q}+\alpha C$, and therefore
$[{\bf u}, {\bf v}]$ belongs to $S(C, {\bf p}, {\bf q}, {\bf x})$.
\hfill{$\Box$}

\medskip
\noindent
{\bf Definition 1.} {\it Let ${\bf p}$ and ${\bf q}$ be fixed distinct points. For
${\bf x}\in H({\bf p}, {\bf q})$ let $\overline{\bf x}$ denote the
middle point of $S(C, {\bf p}, {\bf q}, {\bf x})$. Then we define
$$B(C, {\bf p}, {\bf q})=\left\{ \overline{\bf x}:\ {\bf x}\in H({\bf p},
{\bf q})\right\}$$ and call it a Minkowski bisector of ${\bf pq}$ with
respect to $\| \cdot \|_C$.}

\medskip\noindent
{\bf Remark 1.} {\it In the literature like {\rm \cite{day47}},
{\rm \cite{gru74}}, {\rm \cite{hor00}}, {\rm \cite{jam45}}, {\rm
\cite{man35}}, {\rm \cite{mar04}} and {\rm \cite{woo69}}, the bisector was defined by
$$B'(C, {\bf p}, {\bf q})=\left\{ S(C, {\bf p}, {\bf q}, {\bf x}):\ {\bf x}
\in H({\bf p}, {\bf q})\right\}.$$ Clearly we have
$$B(C, {\bf p}, {\bf q})\subseteq B'(C, {\bf p}, {\bf q}).$$
As one will see from Section $3$ that, for the purpose to study packing and
covering, $B(C, {\bf p}, {\bf q})$ is more natural than $B'(C,
{\bf p}, {\bf q})$.}

Let $f({\bf x})$ be a map from $H({\bf p}, {\bf q})$ to $B(C, {\bf
p}, {\bf q})$ defined by $f({\bf x})=\overline{\bf x}$. Clearly it
is a one-to-one map. However it is less obvious if it is
continuous in general. We note that in some references such as \cite{hor00} and \cite{mar04}
our \lq\lq continuous" is referred as \lq\lq homeomorphic to a hyperplane".
In fact, as one can see from the following lemma and example, the situation is quite complicated.

\medskip\noindent
{\bf Lemma 2.} {\it In $\mathbb{E}^2$, for any given metric $\|\cdot \|_C$
and given distinct points ${\bf p}$ and ${\bf q}$, the map $f({\bf
x})=\overline{{\bf x}}$ from $H({\bf p}, {\bf q})$ to $B(C, {\bf
p}, {\bf q})$ is continuous.}

\medskip\noindent
{\bf Proof.} Without loss of generality, we assume that ${\bf
p}={\bf o}$ and ${\bf q}={\bf e}_2$. Then the set $H({\bf p}, {\bf
q})$ is the $x$-axis. Let $R_1$ denote the set of points ${\bf
x}\in H({\bf p}, {\bf q})$ such that $S(C, {\bf p}, {\bf q}, {\bf
x})$ is a single point and write $R_2=H({\bf p}, {\bf q})\setminus
R_1$.

By repeating partial argument of Lemma 1, it can be deduced that
$R_1$ is closed and if $(x_1,0)\in R_2$, then $(x_2,0)\in R_2$
whenever $|x_2|\ge |x_1|$. In addition, if $x_0=\max \{ x:
(x,0)\in R_1\}$, ${\bf x}_0=(x_0, 0)$ and $f({\bf x}_0)=(x_0,y_0)$, then $f({\bf x})$ (for
$(x,0)\in R_2$ and $x>0$) is the middle point of $[{\bf g}_1,{\bf
g}_2]$, where ${\bf g}_1=(x, {{y_0}\over {x_0}}x)$ and ${\bf
g}_2=(x, 1+{{y_0-1}\over {x_0}}x)$, as shown in Figure 2. In other words $f({\bf x})$ is
an half straight line for $x\geq x_0$. Similarly, one can deal
with the case that $(x,0)\in R_2$ and $x\leq -x_0$.

\begin{figure}[ht]
\centering
\includegraphics[height=4.5cm,width=8.3cm,angle=0]{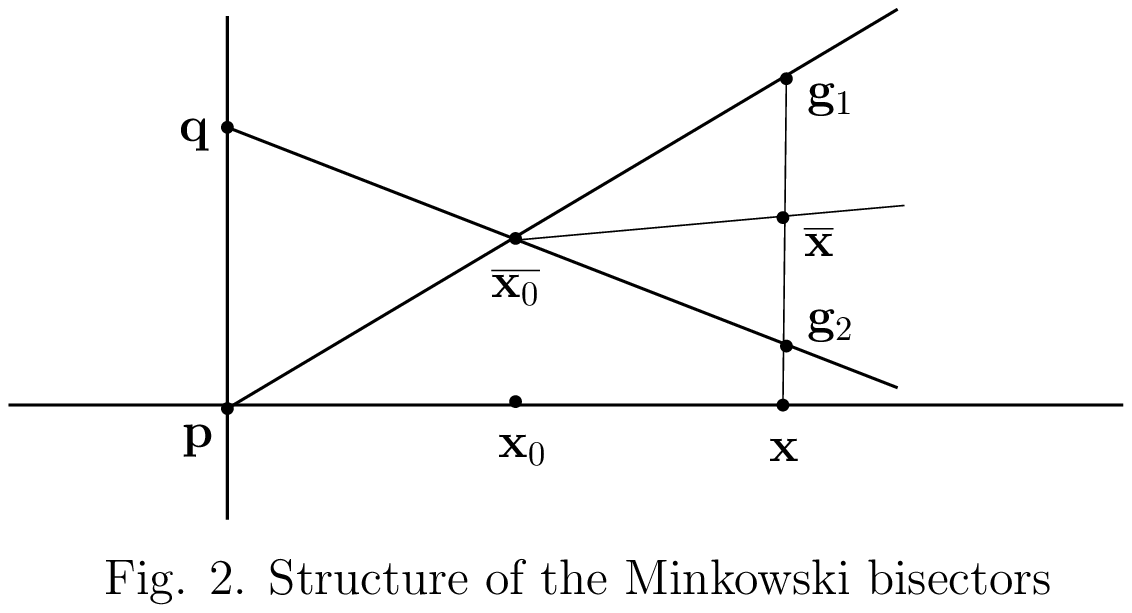}
\end{figure}

Then the continuity of $f({\bf x})$ follows easily, by considering two cases ${\bf x}\in R_1$
and ${\bf x}\in R_2$.
\hfill{$\Box$}

\medskip\noindent
{\bf Example 1.} {\it In $\mathbb{E}^3$ let us define
$$C={\rm conv}\{ {\bf v}, S, -{\bf v}\},$$
where $S=\{ (x,y,0):\ x^2+y^2\le 1\}$ and ${\bf v}=(1,0,1)$. It is
easy to see that $C$ is a centrally symmetric convex body and its
surface has only two vertical segments, they are $[{\bf e}_1, {\bf v}]$ and $[-{\bf e}_1,
-{\bf v}]$. Take ${\bf p}=(0,0,0)$ and ${\bf q}=(0,0,1)$ and let $\mathbb{H}$ denote the set of all planes which containing both ${\bf p}$ and ${\bf q}$. By considering the intersections with planes $H\in \mathbb{H}$, one can deduce that the Minkowski bisector $B(C, {\bf p}, {\bf q})$ with respect to $\| \cdot \|_C$
consists of two parts $L$ and $M$, where $L$ is a
straight line $\{ (x ,0,{1\over 2}(x+1)):\ x \in \mathbb{R}\}$ and
$M$ is a set between two planes $M_1=\{ (x,y,z):\ z=0\}$ and $M_2=\{ (x,y,z):\ z=1\}$, as shown in Figure $3$. Therefore it is easy to see that the map $f({\bf x})=\overline{\bf x}$ is
not continuous at the points $(x, 0,0)$ whenever $|x|>1$. By adding more edges similar to $[{\bf e}_1, {\bf v}]$, one can make the situation much more complicated.}

\begin{figure}[ht]
\centering
\includegraphics[height=6cm,width=8.5cm,angle=0]{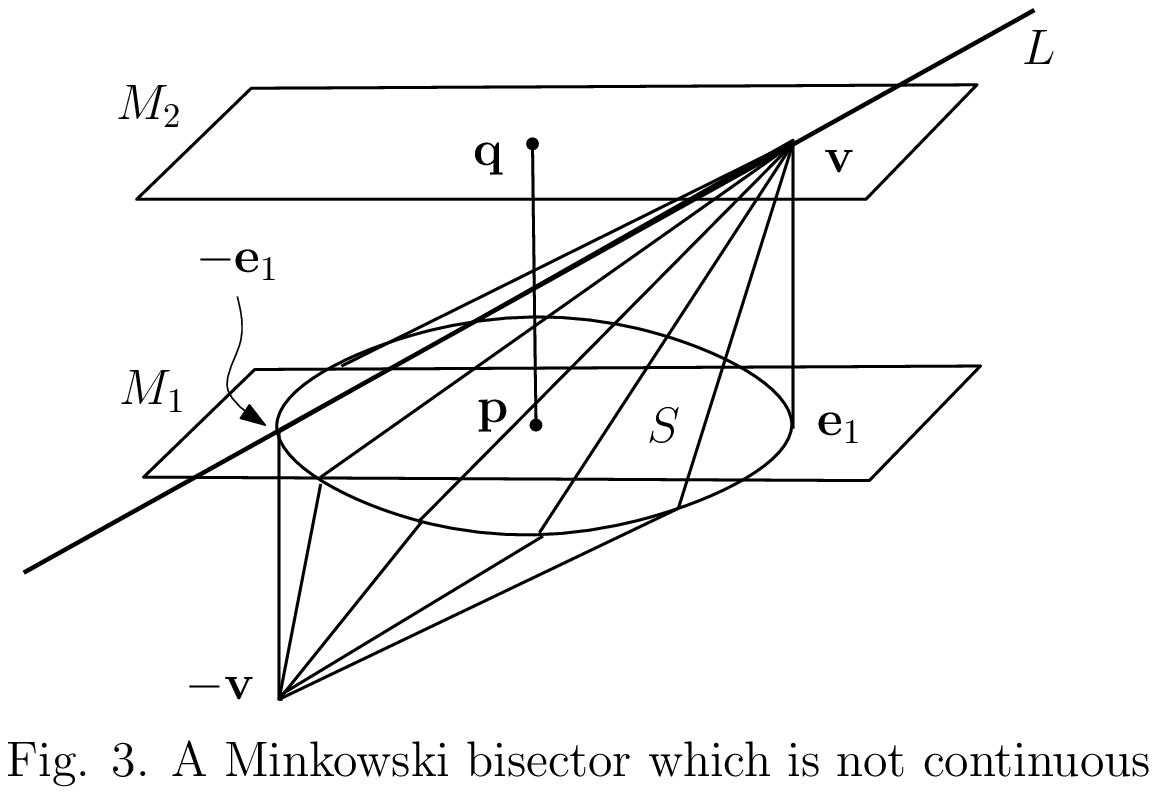}
\end{figure}

\medskip
Let $\mathcal{C}$ denote the set of all $n$-dimensional centrally
symmetric convex bodies centered at the origin ${\bf o}$ and let
$\delta (\cdot, \cdot )$ denote the Hausdorff metric defined on
$\mathcal{C}$. In other words,
$$\delta (C_1, C_2)=\min\{ \gamma :\ C_1\subseteq C_2+\gamma B_n;\
C_2\subseteq C_1+\gamma B_n\}.$$ Let ${\bf p}$ and ${\bf q}$ be
two fixed distinct points, let $C_0$ be a centrally symmetric
convex body, and let $C_1$, $C_2$, $\cdots $ be a sequence of
centrally symmetric convex bodies. It is natural to ask, would
$B(C_i, {\bf p}, {\bf q})$ converges to $B(C_0, {\bf p}, {\bf q})$ if
$$\lim_{i\to\infty }\delta (C_i,C_0)=0\ \! ?$$
Unfortunately, as one will see from the next example, the answer
to this question is negative, even in the plane.

\medskip\noindent
{\bf Example 2.} {\it In the plane we define
$$C_0=\{ (x,y):\ -1\le x\le 1;\ -1\le y\le 1\}$$
and
$$C_i={\rm conv}\left\{ \pm (1,1), \pm (1,0), \pm (1-\mbox{$1\over i$},
-1)\right\}$$ for $i=1, 2, \cdots .$ Then we have
$$\lim_{i\to\infty}\delta (C_i, C_0)=0.$$
Let ${\bf p}=(0,0)$ and ${\bf q}=(0, 2)$, we have
$$B(C_0, {\bf p}, {\bf q})=\{ (x, 1):\ x\in \mathbb{R}\}.$$
However, as shown by Figure $4$, the Minkowski bisector $B(C_i, {\bf p}, {\bf q})$
consists of points $(x,y)$ satisfying
$$ y=\left\{ \begin{array}{ll}
{1\over 2}x+1, & \mbox{if $|x|\ge 2$,}\\
{i\over {i+1}}(x+2), &\mbox{if $-2\le x\le {1\over i}-1$,}\\
1, & \mbox{if $|x|\le 1-{1\over i}$,}\\
{i\over {i+1}}x+{2\over {i+1}}, & \mbox{if $1-{1\over i}\le x\le 2$,}
\end{array}\right.
$$
Therefore the sequence $B(C_i, {\bf p}, {\bf q})$ does not
converge to $B(C_0, {\bf p}, {\bf q}).$}

\begin{figure}[ht]
\centering
\includegraphics[height=5.5cm,width=9cm,angle=0]{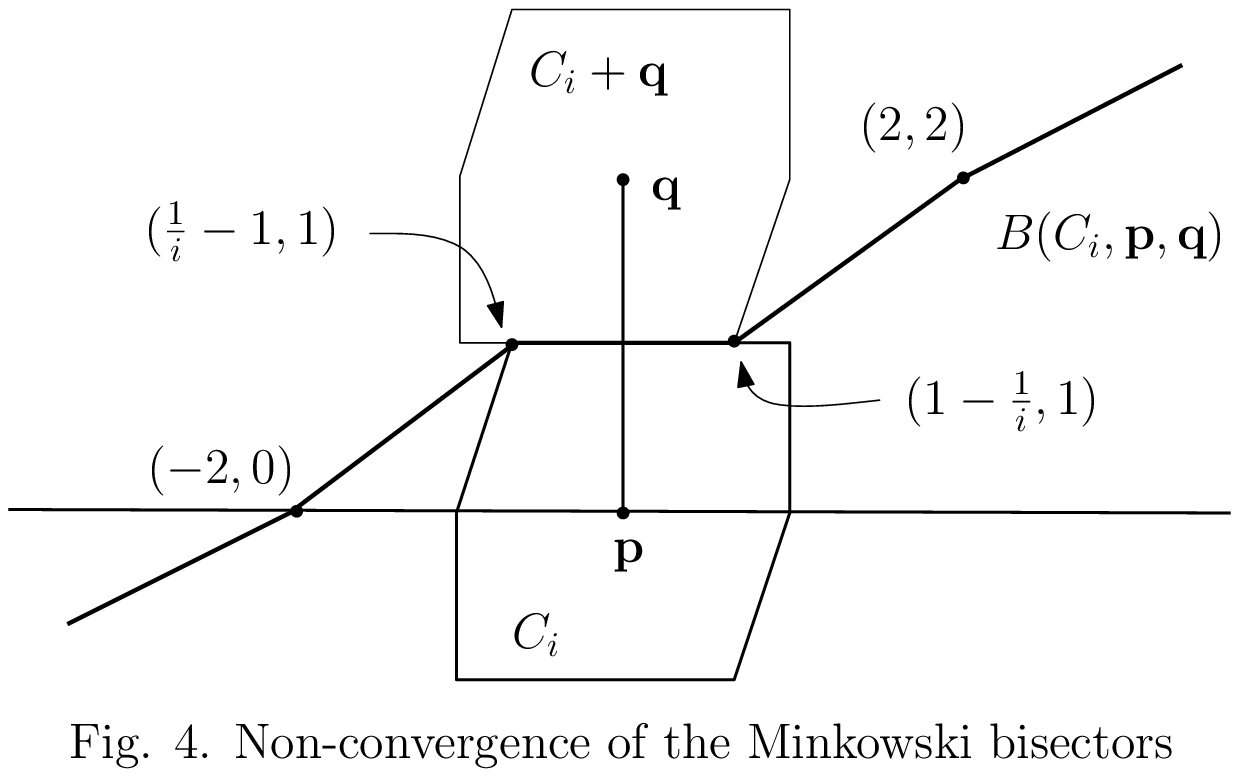}
\end{figure}

\medskip
If $B(C, {\bf p}, {\bf q})$ is a continuous surface, then it
divides the whole space nicely into two half spaces.
Unfortunately, as it was shown by Example 1, the Minkowski bisectors are not
always continuous. Nevertheless, for many important metrics, the
continuity can be guaranteed by the next lemma.

\medskip\noindent
{\bf Lemma 3.} {\it If $C$ is a regular centrally symmetric convex body or a
centrally symmetric polytope, for any pair of distinct points ${\bf p}$ and ${\bf q}$,
the map $f({\bf x})=\overline{{\bf x}}$ from $H({\bf p},{\bf
q})$ to $B(C, {\bf p}, {\bf q})$ is continuous.}

\medskip\noindent
{\bf Proof.} First let us deal with the regular case. By regular
we mean that its surface does not contain any segment. Then, for
any point ${\bf x}\in H({\bf p}, {\bf q})$, $S(C, {\bf p}, {\bf
q}, {\bf x})$ is a single point. If the map is not continuous at
point ${\bf x}_0$, then one can deduce that $S(C, {\bf p}, {\bf
q}, {\bf x}_0)$ is not a single point. The regular case follows.

Now we consider the polytope case. Let $P$ be an $n$-dimensional
centrally symmetric convex polytope, let $P'$ denote its projection
on $H({\bf p}, {\bf q})$, and let $\partial (P')$ denote the relative
boundary of $P'$. Clearly $P'$ is a $(n-1)$-dimensional centrally
symmetric polytope and $\partial (P')$ is a polytope complex of dimension
less than or equal to $n-2$.

Let ${\bf v}$ be a point in $\partial (P')$, let ${\bf v}^*$ denote the point ${\bf
v}+\mu ({\bf q}-{\bf p})\in P$ with maximal $\mu $, let
${\bf v}'$ denote the corresponding point with minimal $\mu $, and write
$$h({\bf v})={{\| {\bf v}^*-{\bf v}'\|}\over {\| {\bf q}-{\bf
p}\|}}.$$ Since $P$ is a polytope, by convexity, it is easy to see
that both $g_1({\bf v})={\bf v}^*$ and $g_2({\bf v})={\bf v}'$
are continuous maps for ${\bf v}\in \partial (P')$, and $h({\bf v})$ is
a continuous function for ${\bf v}\in \partial (P')$.

\begin{figure}[ht]
\centering
\includegraphics[height=5.4cm,width=9cm,angle=0]{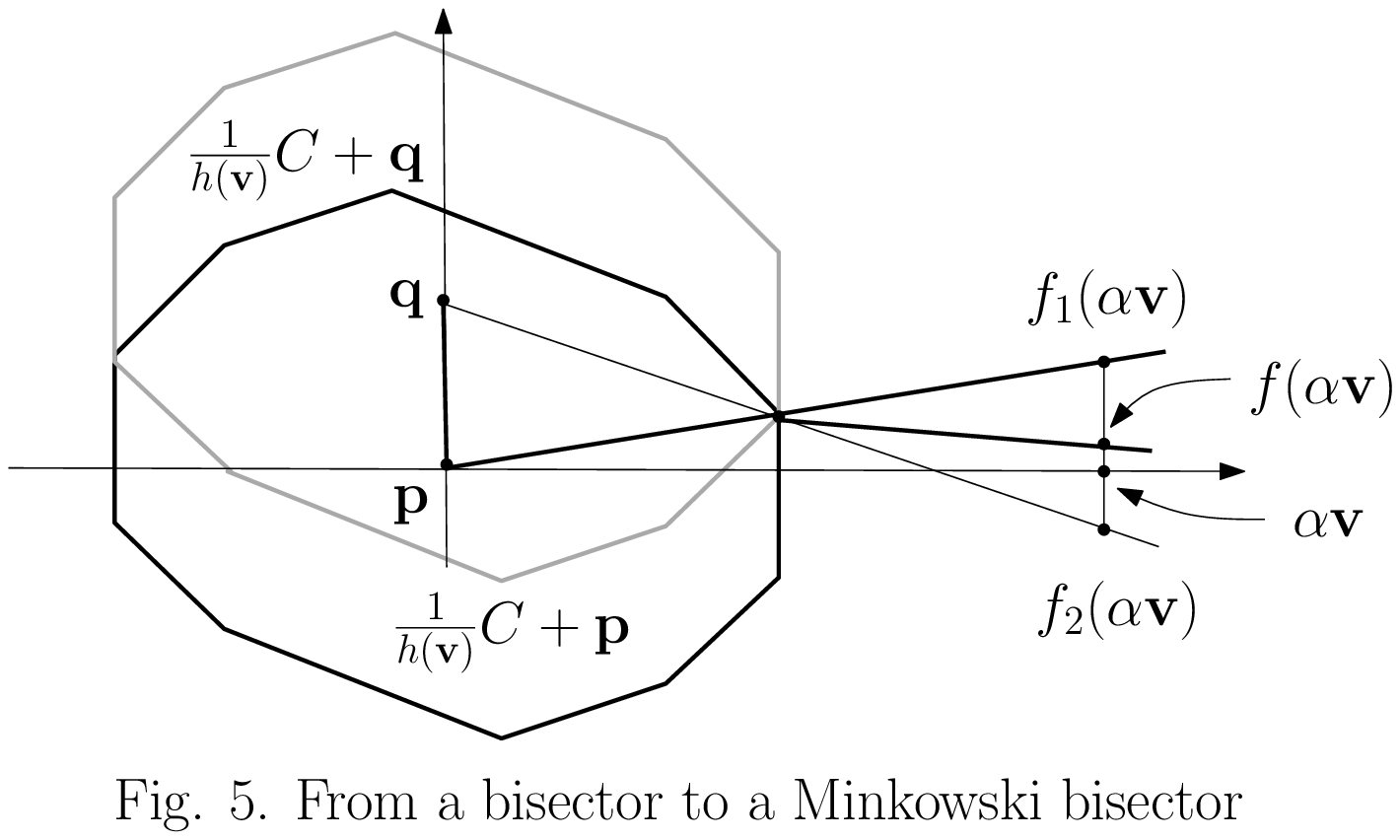}
\end{figure}

Let $Q({\bf p}, {\bf q}, {\bf v})$ denote the two-dimensional hyperplane determined by ${\bf p}$,
${\bf q}$ and ${\bf v}$. By studying the intersections with $Q({\bf p}, {\bf q}, {\bf v})$, it can be
shown that $S(P, {\bf p}, {\bf q}, \alpha {\bf v})$ is a single
point if and only if $\alpha \le 1/h({\bf v})$ (when $h({\bf v})=0$, $\alpha $ can be any number).
For $\alpha {\bf v}\in H({\bf p}, {\bf q})$, let $f_1(\alpha {\bf v})$ denote the
point $\alpha {\bf v}+\mu ({\bf q}-{\bf p})$ in $S(P, {\bf p}, {\bf q}, \alpha {\bf v})$ with the maximal
$\mu$ and let $f_2(\alpha {\bf v})$ denote the corresponding point in $S(P,
{\bf p}, {\bf q}, \alpha {\bf v})$ with the minimal $\mu $. When
$\alpha >1/h({\bf v})$, illustrated by Figure 5,  it can be shown by similar triangles that
$$f_1(\alpha {\bf v})=\alpha {\bf v}^*,$$
$$f_2(\alpha {\bf v})=\alpha {\bf v}^*-(\alpha\ h({\bf v})-1)({\bf q}-{\bf p})$$
and
$$f(\alpha {\bf v})=\alpha {\bf v}^*-\mbox{$1\over 2$}(\alpha\ h({\bf v})-1)({\bf q}-{\bf p}).\eqno (1)$$

Let $H_1({\bf p},{\bf q})$ denote the set of the points ${\bf
x}\in H({\bf p}, {\bf q})$ such that $S(P, {\bf p}, {\bf q}, {\bf
x})$ is not a single point. It can be shown that $H_1({\bf p},{\bf
q})$ is an open set and therefore $H({\bf p}, {\bf q})\setminus
H_1({\bf p}, {\bf q})$ is a closed set. By (1) it follows that
$f({\bf x})$ is continuous in $\overline{H_1({\bf p}, {\bf q})}$.
On the other hand $f({\bf x})$ is continuous in $H({\bf p}, {\bf
q})\setminus H_1({\bf p}, {\bf q})$. Therefore $f({\bf x})$ is a
continuous map from $H({\bf p}, {\bf q})$ to $B(P, {\bf p}, {\bf
q})$. The lemma is proved. \hfill{$\Box$}

\medskip
By the previous proof, it is easy to
see that the Minkowski bisector $B(P, {\bf p}, {\bf q})$ is consists of
planar pieces if the metric is defined by a polytope $P$. Let
$\varphi (P, {\bf p}, {\bf q})$ denote the minimal number $k$ such
that the Minkowski bisector of ${\bf p}$ and ${\bf q}$ with respect to $\|\cdot \|_P$
can be divided into $k$ planar pieces. We have the following
results.

\medskip\noindent
{\bf Theorem 1.} {\it If $P$ is a centrally symmetric polygon with
$m$ vertices. Then, we have} $$\varphi (P, {\bf p}, {\bf q})\le m-1.$$

\medskip\noindent
{\bf Proof.} For convenience, without loss of generality, we take ${\bf p}=(0,0)$ and
${\bf q}=(0,1)$. Let $\lambda$ change from $0$ to infinity, it is easy
to see that the boundaries of $\lambda P+{\bf p}$ and $\lambda
P+{\bf q}$ intersect each other whenever $\lambda \ge \alpha$ for
some $\alpha $.

\begin{figure}[ht]
\centering
\includegraphics[height=5.7cm,width=8.2cm,angle=0]{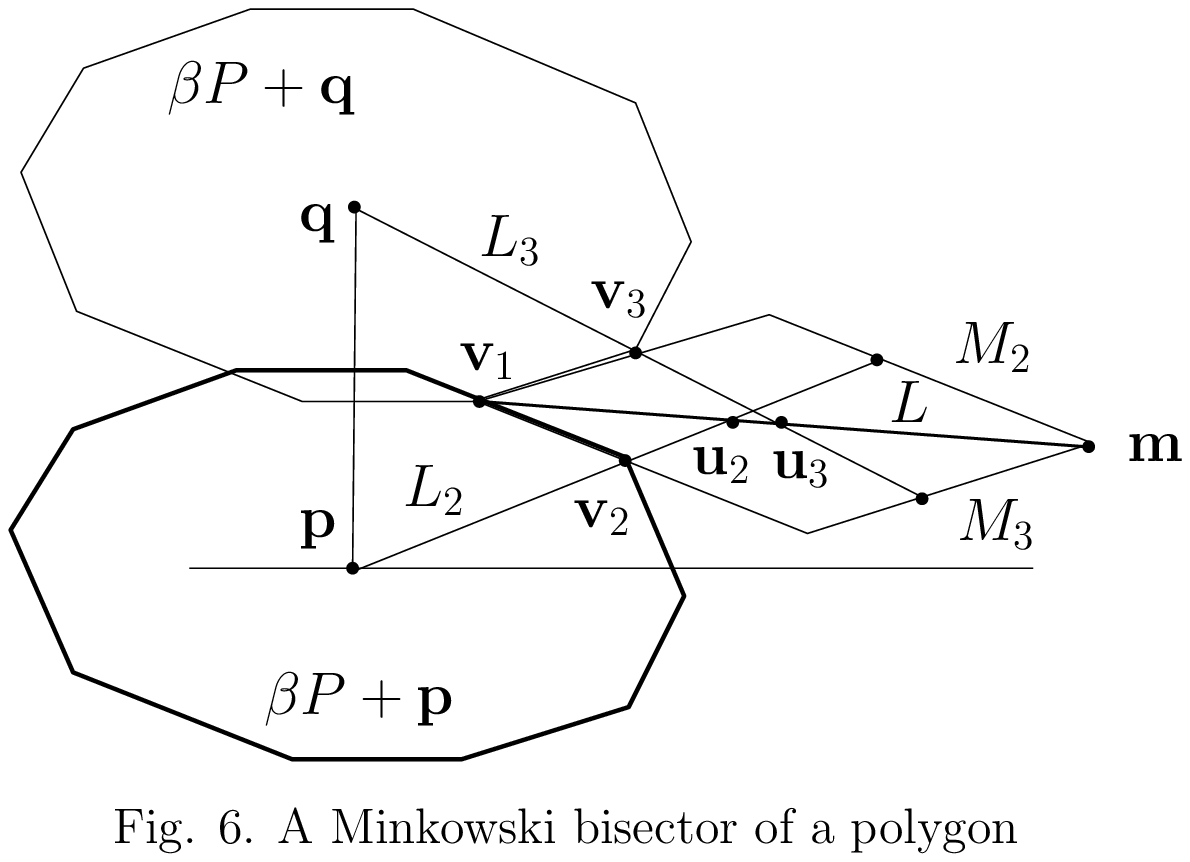}
\end{figure}

For some $\beta >0$, assume that the boundaries of $\beta P+{\bf
p}$ and $\beta P+{\bf q}$ meet at a vertex ${\bf v}_1$ of $\beta
P+{\bf p}$ or $\beta P+{\bf q}$ (see Figure 6), ${\bf v}_1{\bf v}_2$ is an edge of
$\beta P+{\bf p}$, ${\bf v}_1{\bf v}_3$ is an edge of $\beta
P+{\bf q}$. Let $M_2$ denote the line which is parallel to ${\bf v}_1{\bf v}_2$ and pass through
$2{\bf v}_2$, let $M_3$ denote the line which is parallel to ${\bf v}_1{\bf v}_3$ and pass through
$2({\bf v}_3-{\bf q})+{\bf q}$, let ${\bf m}$ be the intersection of $M_2$ and $M_3$, let $L$ denote the line passes ${\bf v}_1$ and ${\bf m}$, let $L_2$ denote the line determined by ${\bf p}$ and ${\bf v}_2$,
and let $L_3$ denote the line determined by ${\bf q}$ and ${\bf
v}_3$. Then $L_2$ intersects $L$ at a point ${\bf u}_2$ and $L_3$
intersects $L$ at a point ${\bf u}_3$. Let ${\bf u}$ denote the
${\bf u}_i$ which is closer to ${\bf v}_1$. By elementary geometry
it is easy to see that the whole segment $[ {\bf u}, {\bf v}_1]$
belongs to $B(P, {\bf p}, {\bf q})$. In addition, for some $\gamma
>0$, the boundaries of $\gamma P+{\bf p}$ and $\gamma P+{\bf q}$
meet at ${\bf u}$ and ${\bf u}$ is a vertex of either $\gamma
P+{\bf p}$ or $\gamma P+{\bf q}$. By repeating this process and dealing with two cases with respect to if $\partial (C)$ contains a segment which is parallel to ${\bf p}{\bf q}$ or not, it can
be shown that
$$\varphi (P, {\bf p}, {\bf q})\le m-1.$$

The theorem is proved. \hfill{$\Box $}

\medskip\noindent
{\bf Remark 2.} {\it It is easy to see that
$$B(\lambda C, {\bf p}, {\bf q})=B(C, {\bf p}, {\bf q})$$
holds for all positive number $\lambda $. Thus, for a fixed polytope $P$, the number $\varphi (P, {\bf p}, {\bf q})$ is determined by the direction of ${\bf p}{\bf q}$. In fact, it follows by the proof of Theorem $1$ that, $$\varphi (P,{\bf p}, {\bf q})=m-1$$ holds for all ${\bf p}$ and ${\bf q}$, except for a finite number of directions ${\bf p}{\bf q}$.}

\bigskip
In higher dimensions the situation is much more complicated. We
have the following upper bound for $\varphi (P, {\bf p}, {\bf q})$.

\medskip\noindent
{\bf Theorem 2.} {\it Let $n\ge 3$ be an integer and let $P$ be an
$n$-dimensional centrally symmetric polytope with $m$ facets. Then
$$\varphi (P, {\bf p}, {\bf q})\le \mbox{$1\over 4$}m^2+\mbox{$1\over {27}$}m^3$$
holds for any pair of distinct points ${\bf p}$ and ${\bf q}$.}

\medskip\noindent {\bf Proof.} Let $F$ denote a facet of $P$, let ${\bf
n}(F)$ denote its outer unit normal, and let $\mathcal{F}$ denote
the set of the $m$ facets of $P$. Then we define
$$\mbox{$\mathcal{F}$}_{-}=\{ F\in \mbox{$\mathcal{F}$}:\ \langle {\bf n}(F),
{\bf q}-{\bf p}\rangle <0\},$$
$$\mbox{$\mathcal{F}$}_{0}=\{ F\in \mbox{$\mathcal{F}$}:\ \langle {\bf n}(F),
{\bf q}-{\bf p}\rangle =0\}\
$$ and
$$\mbox{$\mathcal{F}$}_{+}=\{ F\in \mbox{$\mathcal{F}$}:\ \langle {\bf n}(F),
{\bf q}-{\bf p}\rangle >0\}.$$ It is easy to see, let $|X|$ denote the number of the elements of $X$, that these three sets are pairwise disjoint and satisfying
$$| \mbox{$\mathcal{F}$}_{-}|=| \mbox{$\mathcal{F}$}_{+}|\le \mbox{$1\over 2$}m\eqno (2)$$ and
$$| \mbox{$\mathcal{F}$}_{-}|+ |\mbox{$\mathcal{F}$}_{0}|+ |\mbox{$\mathcal{F}$}_{+}|=m.\eqno (3)$$

Assume that $B(P,{\bf p}, {\bf q})$ can be divided into $\varphi
(P, {\bf p}, {\bf q})$ planar pieces $G_1$, $G_2$, $\cdots $,
$G_{\varphi (P, {\bf p}, {\bf q})}$ and write $$\mathcal{G}=\{
G_i:\ i=1, 2, \cdots , \varphi (P, {\bf p}, {\bf q})\}.$$ For each
$i$ let ${\bf w}_i$ be a relative interior point of $G_i$ such
that there are a relative interior point ${\bf u}_i$ of a facet
which belongs to either $\mathcal{F}_-$ or $\mathcal{F}_0$, a
relative interior point ${\bf v}_i$ of a facet which belongs to
either $\mathcal{F}_0$ or $\mathcal{F}_+$, and a positive number
$\lambda_i$ satisfying
$${\bf w}_i=\lambda_i{\bf u}_i+{\bf q}=\lambda_i{\bf v}_i+{\bf p}.\eqno (4)$$
If ${\bf u}_i\in F\in \mathcal{F}_0$, by looking at Figure 5, it can be deduced that ${\bf v}_i\in F$.

Then $\mathcal{G}$ can be divided into two disjoint subsets
$$ \mathcal{G}_1=\bigg\{ G_i:\ {\bf u}_i\in \bigcup_{F\in
\mathcal{F}_-}F\bigg\}$$ and
$$ \mathcal{G}_2=\bigg\{ G_i:\ {\bf u}_i\in \bigcup_{F\in
\mathcal{F}_0}F\bigg\}.$$ Clearly we have
$$\varphi (P, {\bf p}, {\bf q})=|\mathcal{G}_1| + |\mathcal{G}_2|.\eqno (5)$$ Now we proceed to
estimate $|\mathcal{G}_1|$ and $|\mathcal{G}_2|$, respectively.

If ${\bf u}_i$ and ${\bf u}_j$ belong to one facet $F^*$ in
$\mathcal{F}_-$ and ${\bf v}_i$ and ${\bf v}_j$ belong to one
facet $F'$ in $\mathcal{F}_+$, then one can deduce
that ${\bf w}_i$ and ${\bf w}_j$ should belong to the same planar
piece. In other words, the whole segment $[ {\bf w}_i, {\bf w}_j]$
belongs to $B(P, {\bf p}, {\bf q})$. To see this, when $0\le \theta \le 1$, writing $\lambda =\theta \lambda_i+(1-\theta )\lambda_j$ and $\theta '=\theta \lambda_i /\lambda $, we have
\begin{eqnarray*}
\theta {\bf w}_i+(1-\theta ){\bf w}_j&=& \theta (\lambda_i{\bf u}_i+{\bf q})
+(1-\theta )(\lambda_j{\bf u}_j+{\bf q})\\
&=&\theta \lambda_i{\bf u}_i+(1-\theta )\lambda_j{\bf u}_j+{\bf q}\\
&=&\lambda (\theta'{\bf u}_i+(1-\theta '){\bf u}_j)+{\bf q}
\end{eqnarray*}
and
\begin{eqnarray*}
\theta {\bf w}_i+(1-\theta ){\bf w}_j&=& \theta (\lambda_i{\bf v}_i+{\bf p})
+(1-\theta )(\lambda_j{\bf v}_j+{\bf p})\\
&=&\theta \lambda_i{\bf v}_i+(1-\theta )\lambda_j{\bf v}_j+{\bf p}\\
&=&\lambda (\theta'{\bf v}_i+(1-\theta '){\bf v}_j)+{\bf p},
\end{eqnarray*}
where
$$\theta'{\bf u}_i+(1-\theta '){\bf u}_j\in F^*\quad {\rm and}\quad  \theta'{\bf v}_i+(1-\theta '){\bf v}_j\in F'.$$ Therefore $|\mathcal{G}_1|$ is bounded from above by the number of distinct pairs of facets $\{ F^*, F'\}$ such
that $F^*\in \mathcal{F}_-$ and $F'\in \mathcal{F}_+$. Then by (2) we get
$$| \mathcal{G}_1| \le |\mathcal{F}_-|\cdot |\mathcal{F}_+|\le \mbox{$1\over
4$}m^2.\eqno (6)$$

For $F\in \mathcal{F}_0$ we write
$$ \mathcal{G}(F)=\{ G_i:\ {\bf u}_i\in F\}.$$
Let $\mathcal {F}_+'$ denote the set of the facets $F'$ of $P$ such that $\langle {\bf n}(F'), {\bf q}-{\bf p}\rangle >0$ and $F\cap F'$ is a $(n-2)$-dimensional polytope and let $\mathcal {F}_-'$ denote the set of the facets $F^*$ of $P$ such that $\langle {\bf n}(F^*), {\bf q}-{\bf p}\rangle <0$ and $F\cap F^*$ is a $(n-2)$-dimensional polytope. For ${\bf u}_i$, ${\bf v}_i$ and ${\bf w}_i$ defined by (4), let ${\bf u}^+_i$ denote the
point ${\bf u}_i+\lambda ({\bf q}-{\bf p})\in P$ with the maximal $\lambda $, and let ${\bf u}^-_i$ denote the
point ${\bf u}_i+\lambda ({\bf q}-{\bf p})\in P$ with the minimal $\lambda $. Then, we have
${\bf u}_i^+={\bf v}_i^+$ and ${\bf u}_i^-={\bf v}_i^-$. If both ${\bf u}^+_i$ and ${\bf u}^+_j$ belong to ${\rm int}(F\cap F')$ for a facet $F'\in \mathcal{F}_+'$ and both ${\bf v}^-_i$ and ${\bf v}^-_j$ belong to ${\rm int}(F\cap F^*)$ for a facet $F^*\in \mathcal{F}_-'$, by an argument similar to the previous case it can be deduced that both ${\bf w}_i$ and ${\bf w}_j$ belong to a planar piece of $B(P,{\bf p}, {\bf q})$. Therefore $|\mathcal{G}(F)|$ is bounded from above by the number of distinct pairs of facets $\{ F', F^*\}$ such that $F'\in \mathcal{F}_+'$ and
$F^*\in \mathcal{F}_-'$. Consequently, we get
$$|\mathcal{G}(F)|\le |\mathcal{F}_-'|\cdot |\mathcal{F}_+'|\le |\mathcal{F}_-|\cdot |\mathcal{F}_+|$$
and, by (3),
$$|\mathcal{G}_2| \le |\mathcal{F}_0|\cdot |\mathcal{F}_-|\cdot |\mathcal{F}_+|\le \mbox{${1\over {27}}$}m^3.
\eqno(7)$$

As a conclusion of (5), (6) and (7) we get
$$\varphi (P, {\bf p}, {\bf q})\le \mbox{$1\over 4$}m^2+\mbox{$1\over {27}$}m^3.$$
The theorem is proved.\hfill{$\square $}

\medskip\noindent
{\bf Remark 3.} {\it Although we can not give a proof, the upper
bound in Theorem $2$ seems much too large. It is easy to see that
the direction ${\bf p}{\bf q}$ such that $ \mathcal{F}_0\not=\emptyset $ is a zero measure set on
$\partial (B_n)$. Therefore we have
$$\varphi (P, {\bf p}, {\bf q})\le\mbox{$1\over 4$}m^2,$$
unless the direction ${\bf pq}$ belongs to a zero
measure set of $\partial (B_n)$.}

\medskip
\section{Minkowski Cells, Lattice Packings and Lattice Coverings}

\bigskip
Let $C$ be an $n$-dimensional centrally symmetric convex body, and let ${\bf p}$ and ${\bf q}$ be two distinct points in $\mathbb{E}^n$. We recall that $H({\bf p}, {\bf q})$ is the hyperplane $\{ {\bf x}:\ \langle {\bf x}, {\bf q}-{\bf p}\rangle =0\}$ and $\overline{{\bf x}}$ is the point ${\bf x}+\lambda ({\bf q}-{\bf p})$ in $B(C, {\bf p}, {\bf q})$. For convenience, we define
$$D(C, {\bf p}, {\bf q})=\{\overline{\bf x}+\lambda ({\bf q}-{\bf p}):\ {\bf x}\in H({\bf p}, {\bf q}),\ \lambda \le 0\}.$$
In principle, the structure of $D(C,{\bf p},{\bf q})$ can be very complicated. Nevertheless, we have the following general result.

\medskip\noindent
{\bf Lemma 4.} {\it The set $D(C,{\bf p},{\bf q})$ is a star set with ${\bf p}$ as its origin.}

\medskip
\noindent
{\bf Proof.} For convenience, we take ${\bf p}={\bf o}$. Then we proceed to show that, if
${\bf u}\in D(C, {\bf p},{\bf q})$, then
$\alpha {\bf u}\in D(C,{\bf p}, {\bf q})$ holds for all $\alpha$ with $0\le \alpha \le 1$.

Let $Q({\bf p}, {\bf q}, {\bf u})$ denote the two-dimensional plane determined by ${\bf p}$, ${\bf q}$ and ${\bf u}$. We consider the intersection of $D(C, {\bf p}, {\bf q})$ with $Q({\bf p}, {\bf q}, {\bf u})$. Clearly, ${\bf u}\in D(C, {\bf p},{\bf q})$ implies
$$\| {\bf p}, {\bf u}\|_C\le \| {\bf q}, {\bf u}\|_C.\eqno(8)$$

\begin{figure}[ht]
\centering
\includegraphics[height=5.3cm,width=7cm,angle=0]{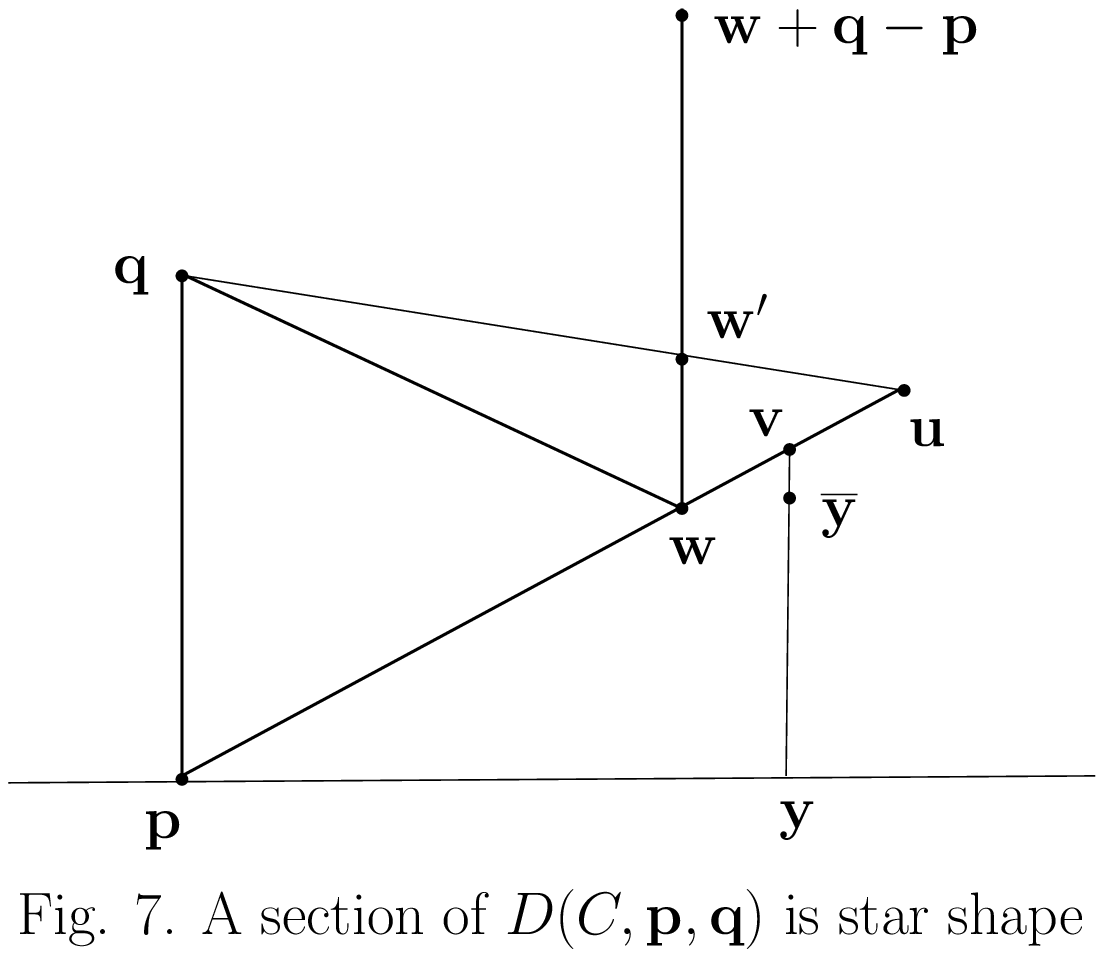}
\end{figure}

On the contrary, if
$${\bf v}=\alpha {\bf u}\not\in D(C,{\bf p},{\bf q})$$
holds for some $\alpha <1$, we proceed to deduce a contradiction. Assume that
$${\bf v}={\bf y}+\lambda ({\bf q}-{\bf p})$$
with suitable ${\bf y}\in H({\bf p},{\bf q})$ and $\lambda \in \mathbb{R}.$ Then the point $\overline{\bf y}$ corresponding to ${\bf y}$ on the Minkowski bisector is below ${\bf v}$, as shown in Figure 7. By Lemma 2, there is a suitable ${\bf w}$ between ${\bf p}$ and ${\bf v}$ such that ${\bf w}\in B(C,{\bf p},{\bf q})$. Therefore,
$${\bf w}\in \partial (\lambda C+{\bf p})\cap \partial (\lambda C+{\bf q})$$
holds for some suitable $\lambda $. Then we have
$${\bf w}+{\bf q}-{\bf p}\in \partial (\lambda C+{\bf q}).$$

If the whole segment $[{\bf w}, {\bf w}+{\bf q}-{\bf p}]$ belongs to $\partial (\lambda C+{\bf q})$, then one can deduce that the whole segment $[{\bf w}, {\bf u}]$ belongs to $D(C,{\bf p}, {\bf q})$ and thus ${\bf v}\in D(C,{\bf p},{\bf q})$, which contradicts the assumption. If
$$[{\bf w}, {\bf w}+{\bf q}-{\bf p}]\not\subset\partial (\lambda C+{\bf q}),$$
then by convexity it can be deduced that
$$\| {\bf q}, {\bf u}\|_C<{{\| {\bf q}, {\bf u}\|}\over {\| {\bf q}, {\bf w}'\|}}\cdot \lambda
={{\| {\bf p}, {\bf u}\|}\over {\| {\bf p}, {\bf w}\|}}\cdot \lambda =\| {\bf p}, {\bf u}\|_C,$$
which contradicts (8).

As a conclusion, $\alpha {\bf u}\in D(C,{\bf p},{\bf q})$ holds for all $0\le \alpha \le 1.$ Lemma 4 is proved. \hfill{$\square$}

\medskip
\noindent
{\bf Remark 4.} {\it For the $C$, ${\bf p}$ and ${\bf q}$ defined in Example $1$, the region $D(C, {\bf p}, {\bf q})$ is not closed. However, when $B(C, {\bf p}, {\bf q})$ is continuous, the region $D(C, {\bf p}, {\bf q})$ is closed.}

\medskip
\noindent
{\bf Definition 2.} {\it Let $C$ be a centrally symmetric convex body, let $\Lambda $ be a lattice and define
$$M(C, \Lambda )=\bigcap_{{\bf v}\in \Lambda \setminus \{ {\bf o}\}}D(C, {\bf o}, {\bf v}).$$
We call $M(C, \Lambda )$ a Minkowski cell of $\Lambda $ with respect to the metric $\| \cdot\|_C$. }

\medskip
\noindent
{\bf Remark 5.} {\it For fixed $C$, ${\bf p}$ and ${\bf q}$, it can be verified that both
$$B(\lambda C, {\bf p}, {\bf q})=B(C, {\bf p}, {\bf q})$$
and
$$D(\lambda C, {\bf p}, {\bf q})=D(C, {\bf p}, {\bf q})$$
hold for all positive numbers $\lambda $. If $\Lambda $ is a lattice and $\tau $ is a non-singular linear transformation from $\mathbb{E}^n$ to $\mathbb{E}^n$. Then we have
$$B(\tau (C), \tau ({\bf p}), \tau ({\bf q}))=\tau (B(C, {\bf p}, {\bf q})),$$
$$D(\tau (C), \tau ({\bf p}), \tau ({\bf q}))=\tau (D(C, {\bf p}, {\bf q}))$$
and
$$M(\tau (C), \tau (\Lambda ))=\tau (M(C, \Lambda )).$$
By Lemma $4$ it follows that the Minkowski cells are centrally symmetric star sets centered at the origin. In the Euclidean case, $C=B_n$, all Minkowski bisectors are hyperplanes, and all Minkowski cells are parallelohedra.}

\medskip
\noindent
{\bf Lemma 5.} {\it Let $C$ be an $n$-dimensional centrally symmetric convex body, let $\Lambda $ be an $n$-dimensional lattice, and let $\gamma $ denote the smallest positive number such that $\gamma C+\Lambda $ is a covering of $\mathbb{E}^n$. Then}
$$M(C,\Lambda )\subseteq \gamma C.$$

\smallskip\noindent
{\bf Proof.} If, on the contrary, there is a point ${\bf x}\in M(C, \Lambda )$ such that
$\| {\bf o}, {\bf x}\|_C>\gamma .$ Since $\gamma C+\Lambda =\mathbb{E}^n$, there is a lattice point ${\bf v}\in \Lambda$ such that $\| {\bf v}, {\bf x}\|_C\le \gamma .$ Thus, we have
$${\bf x}\not\in D(C, {\bf o}, {\bf v})$$
and hence
$${\bf x}\not\in M(C, \Lambda ).$$
The lemma is proved. \hfill{$\square$}

\medskip
\noindent
{\bf Theorem 3.} {\it Let $C$ be a two-dimensional centrally symmetric convex domain and let $\Lambda $ be a two-dimensional lattice. Then the region $M(C, \Lambda )$ is a centrally symmetric star domain and $M(C,\Lambda )+\Lambda $ is a tiling.}

\medskip
\noindent
{\bf Proof.} First, by Definition 2, Lemma 4 and Remark 4, it follows that $M(C, \Lambda )$ is a centrally symmetric star domain. In other words, it is a centrally symmetric compact star set.

Now, we claim that {\it
$$( {\rm int} (M(C, \Lambda ))+{\bf v}_i)\cap  ({\rm int} (M(C, \Lambda ))+{\bf v}_j)=\emptyset $$
holds for all distinct lattice points ${\bf v}_i$ and ${\bf v}_j$.}

If, on the contrary, without loss of generality there are a point ${\bf x}$, a positive number $\epsilon $
and a lattice point ${\bf v}$ such that
$$\epsilon C+{\bf x}\subseteq  {\rm int} (M(C, \Lambda ))\cap  ({\rm int} (M(C, \Lambda ))+{\bf v}).
\eqno (9)$$
Then, we observe the Minkowski bisector $B(C, {\bf o}, {\bf v})$. If ${\bf x}\in B(C, {\bf o}, {\bf v})$,
then we have
$${\bf x}+\epsilon \overline{{\bf v}}\not\in D(C, {\bf o}, {\bf v})$$
and
$${\bf x}+\epsilon \overline{{\bf v}}\not\in M(C, \Lambda ),$$
which contradicts to (9), where $\overline{{\bf v}}$ is the boundary point of $C$ in the direction of ${\bf v}$ as defined in Section 1. If ${\bf x}\not\in B(C, {\bf o}, {\bf v})$, since $B(C, {\bf o}, {\bf v})$ divides $\mathbb{E}^2$ into two separated parts which contains ${\rm int}(M(C, \Lambda ))$ and ${\rm int}(M(C, \Lambda ))+{\bf v}$, respectively. Therefore, $\epsilon C+{\bf x}$ can't belong to both ${\rm int}(M(C, \Lambda ))$ and ${\rm int}(M(C, \Lambda ))+{\bf v}$ simultaneously when $\epsilon $ is sufficiently small, which contradicts to (9) as well.

Next, we claim that {\it for every point ${\bf x}\in \mathbb{E}^2$ there is a lattice point ${\bf v}$ such that}
$${\bf x}\in M(C, \Lambda )+{\bf v}.\eqno (10)$$

Without loss of generality, since $\Lambda $ is periodic, we assume that ${\bf x}\in \gamma C$. Clearly, we have
$$\gamma C\subseteq D(C, {\bf o}, {\bf q})$$
whenever $\| {\bf o}, {\bf q}\|_C>2\gamma.$ Therefore, there are at most $|2\gamma C\cap \Lambda |$ Minkowski bisectors which can effect $M(C, \Lambda )$. Consequently, we assume further that
$${\bf x}\in {\rm int}(\gamma C)\setminus \bigcup_{{\bf v}_i, {\bf v}_j\in 2\gamma C\cap \Lambda }B(C, {\bf v}_i, {\bf v}_j).$$

If, for any positive $\epsilon $ there is a point ${\bf x}'$ such that ${\bf x}'\in \epsilon C+{\bf x}$ and
$$\| {\bf o}, {\bf x}'\|_C< \| {\bf v}, {\bf x}'\|_C$$
holds for all ${\bf v}\in 2\gamma C\cap (\Lambda\setminus \{ {\bf o}\}).$ Then, since $M(C, \Lambda )$ is
compact, we have ${\bf x}\in M(C, \Lambda )$. If, there are a positive number $\epsilon '$ and a subset $W$ of $2\gamma C\cap \Lambda$ with $|W|\ge 2$ such that
$$\| {\bf x}', {\bf w}_i\|_C=\min \{ \|{\bf x}', {\bf v}\|_C :\ {\bf v}\in 2\gamma C\cap \Lambda \}$$
holds for all ${\bf x}'\in \epsilon 'C+{\bf x}$ and ${\bf w}_i\in W$. Then, $W\subset \partial (\lambda C)+{\bf x}$ holds for some suitable positive number $\lambda $. By Figure 8 it can be shown that, if ${\bf w}_i$ and ${\bf w}_j$
are two distinct points in $W$, then the whole segment $[{\bf w}_i, {\bf w}_j]$ belongs to $\partial (\lambda C)+{\bf x}$ and thus all the points of $W$ are colinear.

\begin{figure}[ht]
\centering
\includegraphics[height=5.3cm,width=7cm,angle=0]{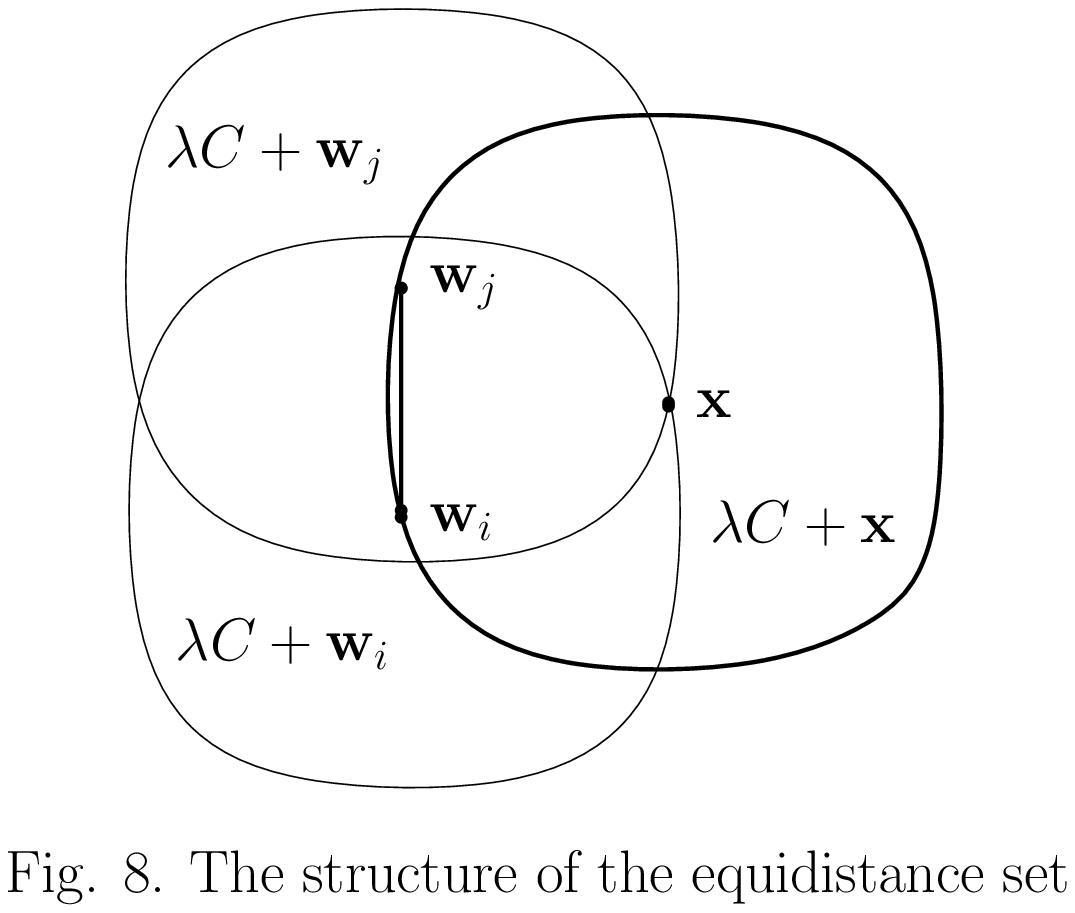}
\end{figure}

Assume that ${\bf w}_1$, ${\bf w}_2$, $\cdots $, ${\bf w}_k$ are the points of $W$ that successively on a line
$H$, as shown in Figure 9. Since $W\subset \Lambda $, we have
$${\bf w}_{i+1}-{\bf w}_i={\bf w}_2-{\bf w}_1.$$
Let $\Lambda '$ denote the affine lattice generated by $W$ in $H$. In other words,
$$\Lambda'= \left\{ {\bf w}_1+z ({\bf w}_2-{\bf w}_1):\ z\in \mathbb{Z} \right\}.$$
Then, $M(C, \Lambda')+\Lambda'$ is a tiling of $\mathbb{E}^2$.
Therefore, as illustrated by Figure 9, there is a ${\bf w}_i\in W$ such that
$${\bf x}\in M(C, \Lambda')+{\bf w}_i$$
and consequently
$${\bf x}\in M(C, \Lambda)+{\bf w}_i.$$
\begin{figure}[ht]
\centering
\includegraphics[height=5.5cm,width=8.5cm,angle=0]{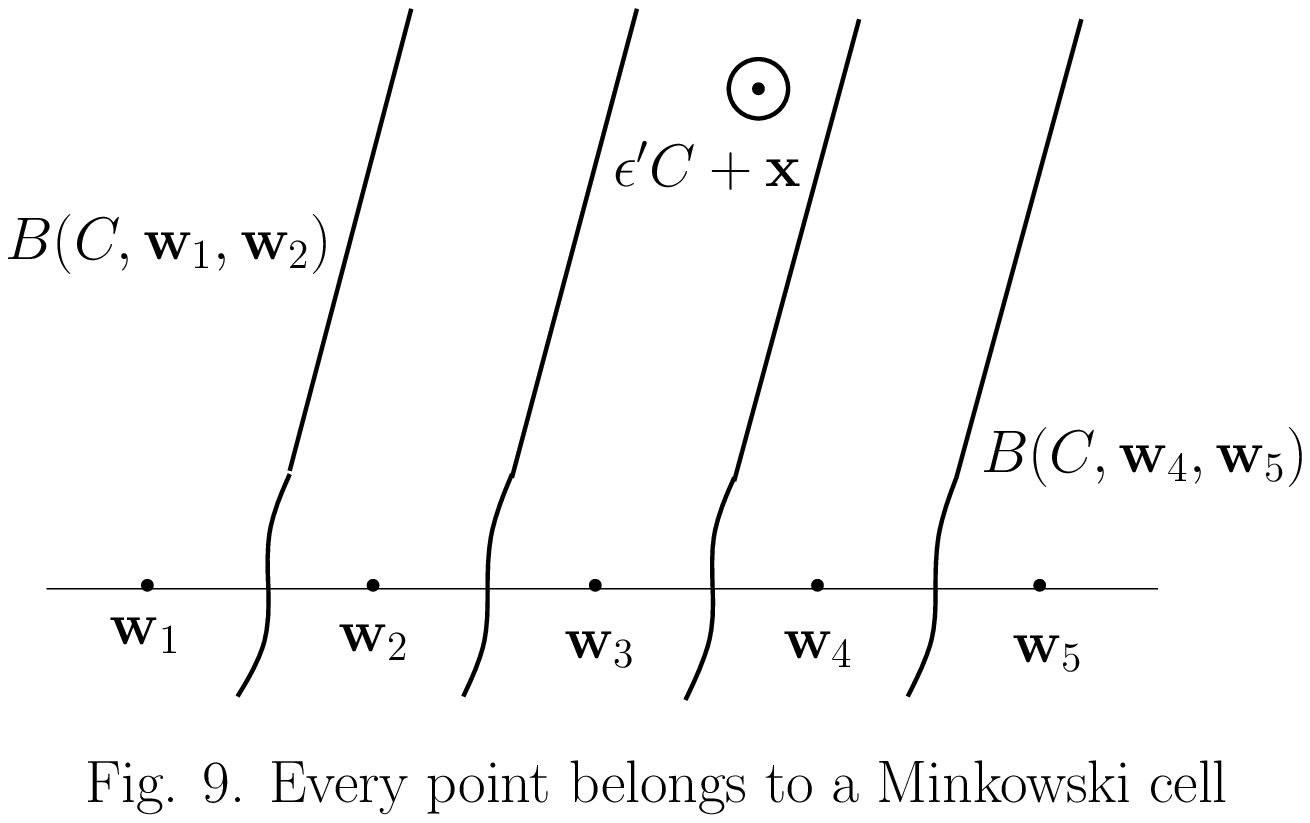}
\end{figure}

As a conclusion, $M(C, \Lambda )+\Lambda $ is a tiling of $\mathbb{E}^2$. Theorem 3 is proved.  \hfill{$\square$}

\medskip
\noindent
{\bf Theorem 4.} {\it Let $C$ be an $n$-dimensional centrally symmetric convex body, let $\Lambda $ be an $n$-dimensional lattice, and let $\gamma $ denote the smallest positive number such that $\gamma C+\Lambda $ is a covering of $\mathbb{E}^n$. If $\partial (C)$ has no segment in the directions of $\{{\bf v}: {\bf v}\in 2\gamma C\cap (\Lambda \setminus \{ {\bf o}\})\}$, then the region $M(C, \Lambda )$ is a centrally symmetric star body and $M(C,\Lambda )+\Lambda $ is a tiling.}

\medskip
This result can be proved just like the Euclidean case. Let $\delta^*(C)$ denote the density of the densest lattice packing of $C$ and let $\theta^*(C)$ denote the density of the thinnest lattice covering of $\mathbb{E}^n$ by $C$. Theorem 4 has the following corollary.

\medskip\noindent
{\bf Corollary 1.} {\it Let $P$ be an $n$-dimensional centrally symmetric polytope and let $\mathcal{M}$ denote the set of all Minkowski cells $M(P, \Lambda )$ contained in $P$. Then, we have}
$$\theta^*(P)=\inf_{M\in \mathcal{M}}{{{\rm vol}(P)}\over {{\rm vol}(M)}}.$$

\medskip\noindent
{\bf Remark 6.} {\it When $C=B_3$, there are only five types of Minkowski cells $($parallelohedra$)$. Therefore, one can determine the values of $\delta^*(B_3)$ and $\theta^*(B_3)$ by studying a unit ball inscribed in parallelohedra or parallelohedra inscribed in a unit ball. For other particular nontrivial centrally symmetric convex bodies, for example the octahedron, to enumerate their Minkowski cells seems challenging and interesting.}

\medskip
It was proved by D.G. Ewald, D.G. Larman and C.A. Rogers \cite{ewal70} that the line segment directions on the surface of an $n$-dimensional convex body is a very small subset of $\partial (B_n)$. Therefore, Theorem $4$ and Corollary 1 seem can be improved further. We end this article by three open problems as following.

\medskip\noindent
{\bf Problem 1.} {\it Let $P$ be an $n$-dimensional centrally symmetric polytope. Enumerate the different types $($geometric or combinatorial$)$ of the Minkowski cells $M(P, \Lambda )$ for all lattices $\Lambda $.}

\medskip\noindent
{\bf Problem 2.} {\it Let $C$ be an $n$-dimensional centrally symmetric convex body such that all its Minkowski bisectors are continuous. Is $M(C, \Lambda )+\Lambda$ always a tiling of $\mathbb{E}^n$ for all lattice $\Lambda $?}

\medskip\noindent
{\bf Problem 3.} {\it Whenever $C+\Lambda$ is a lattice covering of $\mathbb{E}^n$, $n\ge 3$, is there always a parallelohedron $P$ satisfying both $P\subseteq C$ and $P+\Lambda $ is a tiling of $\mathbb{E}^n$?}

\vspace{0.5cm}\noindent
{\bf Acknowledgements.} This work is supported by 973 Programs 2013CB834201 and 2011CB302401, the National Science Foundation of China (No.11071003), and the Chang Jiang Scholars Program of China. For some useful comments, I am grateful to Prof. Senlin Wu.

\bibliographystyle{amsplain}

\vspace{0.5cm}
Chuanming Zong

School of Mathematical Sciences

Peking University

Beijing 100871, China

cmzong@math.pku.edu.cn

\end{document}